\input amstex\documentstyle{amsppt}  
\pagewidth{12.5cm}\pageheight{19cm}\magnification\magstep1
\topmatter
\title From unipotent classes to conjugacy classes in the Weyl group\endtitle
\author G. Lusztig\endauthor
\address{Department of Mathematics, M.I.T., Cambridge, MA 02139}\endaddress
\thanks{Supported in part by the National Science Foundation}\endthanks
\endtopmatter   
\document

\define\uuG{\un{\un G}}

\define\uWW{\un\WW}

\define\si{\sim}

\define\part{\partial}

\define\m{\mapsto}
\define\do{\dots}

\define\T{\times}
\define\ti{\tilde}
\define\nl{\newline}
\redefine\i{^{-1}}

\define\un{\underline}

\define\g{\gamma}

\define\e{\epsilon}

\define\io{\iota}

\define\p{\pi}

\define\ps{\psi}
\define\r{\rho}
\define\s{\sigma}
\redefine\t{\tau}

\define\k{\kappa}

\define\z{\zeta}

\define\Ph{\Phi}
\define\Ps{\Psi}

\define\kk{\bold k}

\define\nn{\bold n}

\define\CC{\bold C}

\define\NN{\bold N}

\define\WW{\bold W}
\define\ZZ{\bold Z}

\define\ce{\Cal E}

\define\cp{\Cal P}
\define\cq{\Cal Q}
\define\car{\Cal R}
\define\cs{\Cal S}
\define\ct{\Cal T}

\define\te{\ti e}

\define\tA{\ti A}

\define\tcp{\ti{\cp}}

\define\CA{Ca}
\define\KL{KL}
\define\UNI{L1}
\define\WE{L2}
\define\SPA{Sp}

\head Introduction\endhead
\subhead 0.1\endsubhead
Let $G$ be a connected reductive algebraic group over an algebraically closed field $\kk$ of characteristic 
$p\ge0$. Let $\uuG$ be the set of unipotent conjugacy classes in $G$. Let $\uWW$ be the set of conjugacy classes 
in the Weyl group $\WW$ of $G$. Let $\Ph:\uWW@>>>\uuG$ be the (surjective) map defined in \cite{\WE}. For 
$C\in\uWW$ we denote by $m_C$ the dimension of the fixed point space of $w:V@>>>V$ where $w\in C$ and $V$ is the 
reflection representation of the Coxeter group $\WW$. The following result is an attempt to construct a
one sided inverse for $\Ph$.

\proclaim{Theorem 0.2} Assume that either $p$ is not a bad prime for $G$ or that $G$ is a simple exceptional 
group. Then for any $\g\in\uuG$ the function $\Ph\i(\g)@>>>\NN$, $C\m m_C$ reaches its minimum at a unique 
element $C_0\in\Ph\i(\g)$. Thus we have a well defined map $\Ps:\uuG@>>>\uWW$, $\g\m C_0$ such that 
$\Ph\Ps:\uuG@>>>\uuG$ is the identity map.
\endproclaim
It is likely that (i) the theorem holds without any assumption and (ii) when $p=0$, the map $\Ps$ coincides with 
the map defined in \cite{\KL, Section 9} (note that the last map has not been computed explicitly in all cases). 
It is enough to prove the theorem in the case where $G$ is almost 
simple; moreover in that case it is enough to consider one group in each isogeny class. When $G$ has type $A$, 
the theorem is trivial since $\Ph$ is a bijection. For $G$ of type $B,C,D$ the proof is given in Section 1. For 
$G$ of exceptional type the proof is given in Section 2.

\head 1. Type $B,C,D$\endhead
\subhead 1.1\endsubhead
Let $\cp^1$ be the set of sequences $c_*=(c_1\ge c_2\ge\do\ge c_m)$ in $\ZZ_{>0}$. For $c_*\in\cp^1$ we set
$|c_*|=c_1+c_2+\do+c_\t$, $\t_{c_*}=m$. For $N\in\NN$ let $\cp^1_N=\{c_*\in\cp^1;|c_*|=N\}$.

Let $\ti\cp=\{p_*\in\cp^1;\t_{p_*}\in2\NN;p_1=p_2,p_3=p_4,\do\}$; for $m\in\NN$ we set
$\ti\cp_{2m}=\ti\cp\cap\cp^1_{2m}$.

Let $\ct$ be the set of $c_*=(c_1\ge c_2\ge\do\ge c_\t)\in\cp^1$ such that for any odd $j$, 
$|\{k\in[1,\t];c_k=j\}|$ is even; for $m\in\NN$ let $\ct_{2m}=\ct\cap\cp^1_{2m}$.

Let $\cs$ be the set of $r_*=(r_1\ge r_2\ge\do\ge r_\s)\in\cp^1$ such that $r_k\in2\ZZ_{>0}$ for all $r$; for 
$m\in\NN$ let $\cs_{2m}=\cs\cap\cp^1_{2m}$.

In this subsection we fix $\nn\in2\NN$.
Let $A'_\nn$ be the set of all pairs $(r_*,p_*)\in\cs\T\ti\cp$ such that $|r_*|+|p_*|=\nn$. 
We define $\io:A'_\nn@>>>\ct_\nn$ by $(r_*,p_*)\m c_*$ where the multiset of entries of $c_*$ is the union of
the multiset of entries of $r_*$ with the multiset of entries of $p_*$.
Let $c_*=(c_1\ge c_2\ge\do\ge c_\t)\in\ct_\nn$. We associate to 
$c_*$ an element $(r_*,p_*)\in\cs\T\ti\cp$ by specifying the number of times $M_e$ (resp. $N_e$) that an integer 
$e\ge1$ appears in $r_*$ (resp. $p_*$). Let $Q_e$ be the number of times that $e$ appears as an entry of $c_*$.
If $e\in2\NN+1$ then $M_e=0$, $N_e=Q_e$. If $e\in2\NN+2$ then $M_e=Q_e,N_e=0$.
Clearly, $c_*\m(r_*,p_*)$ is a well defined map $\io':\ct_\nn@>>>A'_\nn$; moreover, 
$\io\io':\ct_\nn@>>>\ct_\nn$ is the identity map. 

We preserve the notation for $c_*,r_*,p_*$ as above (so that $(r_*,p_*)\in\io\i(c_*)$) and we assume that
$(r'_*,p'_*)\in\io\i(c_*)$. 
Let $M'_e$ (resp. $N'_e$) be the number of times that an integer $e\ge1$ appears in 
$r'_*$ (resp. $p'_*$). Note that $M'_e+N'_e=M_e+N_e$.
If $e\in2\NN+1$ then $M'_e=0$ hence $N'_e=Q_e$. If $e\in2\NN+2$ then $N'_e\ge0$.
We see that in all cases we have $N'_e\ge N_e$. It follows that $\sum_eN'_e\ge\sum_eN_e$ (and the equality 
implies that $N'_e=N_e$ for all $e$ hence $(r'_*,p'_*)=(r_*,p_*)$). We see that 

(a) {\it for any $c_*\in\ct_\nn$ there is exactly one element $(r_*,p_*)\in\io\i(c_*)$ such that the number of
entries of $p_*$ is minimal (that element is $\io'(c_*)$).}

\subhead 1.2\endsubhead
Let $\cp^0=\{p_*\in\cp^1;\t_{p_*}=\text{ even}\}$; for $n\in\NN$ we set $\cp_n^0=\cp^0\cap\cp^1_n$.

Let $\cq$ be the set of all $c_*=(c_1\ge c_2\ge\do\ge c_\t)\in\cp^1$ such that for any even $j$,
$|\{k\in[1,\t];c_k=j\}|$ is even; for $\nn\in\NN$ let $\cq_\nn=\cq\cap\cp^1_\nn$.

Let $\car$ be the set of all $r_*=(r_1\ge r_2\ge\do\ge r_\t)\in\cq$ such that the following conditions are
satisfied. Let $J_{r_*}=\{k\in[1,\t];r_k\text{ is odd}\}$. We write the multiset $\{r_k;k\in J_{r_*}\}$ as a 
sequence $r^1\ge r^2\ge\do\ge r^s$. (We have necessarily $\t=s\mod2$.) We require that:

-if $\t\ne0$ then $1\in J_{r_*}$;

-if $\t\ne0$ is even then $\t\in J_{r_*}$;

-if $u\in[1,s-1]$ is odd then $r^u>r^{u+1}$;

-if $u\in[1,s-1]$ is even then there is no $k'\in[1,\t]$ such that $r^u>r_{k'}>r^{u+1}$.
\nl
For $\nn\in\NN$ we set $\car_\nn=\car\cap\cq_\nn$.

We now fix $\nn\in\NN$. Define $\k\in\{0,1\}$ and $n\in\NN$ by $\nn=2n+\k$. Define $\Xi:\cp_n^\k@>>>\car_\nn$ by
$$(p_1\ge p_2\ge\do\ge p_\s)\m(2p_1+\ps(1)\ge2p_2+\ps(2)\ge\do\ge2p_\s+\ps(\s))$$
if $\s+\k$ is even,
$$(p_1\ge p_2\ge\do\ge p_\s)\m(2p_1+\ps(1)\ge2p_2+\ps(2)\ge\do\ge2p_\s+\ps(\s)\ge1)$$
if $\s+\k$ is odd, where $\ps:[1,\s]@>>>\{-1,0,1\}$ is as follows:

if $t\in[1,\s]$ is odd and $p_t<p_x$ for any $x\in[1,t-1]$ then $\ps(t)=1$;

if $t\in[1,\s]$ is even and $p_x<p_t$ for any $x\in[t+1,\s]$, then $\ps(t)=-1$;

for all other $t\in[1,\s]$ we have $\ps(t)=0$.
\nl
Now $\Xi$ is a bijection with inverse map $\Xi':\car_\nn@>>>\cp_n^\k$ given by
$$(r_1\ge r_2\ge\do\ge r_\t)\m((r_1+\z(1))/2\ge(r_2+\z(2))/2\ge\do\ge(r_\t+\z(\t))/2))$$
if $r_\t>\k$,
$$(r_1\ge r_2\ge\do\ge r_\t)\m((r_1+\z(1))/2\ge(r_2+\z(2))/2\ge\do\ge(r_{\t-1}+\z(\t-1))/2))$$
if $r_\t=\k$, 
where $\z:[1,\t]@>>>\{-1,0,1\}$ is given by $\z(k)=(-1)^k(1-(-1)^{r_k})/2$. (Thus, $\z(k)=(-1)^k$ if $r_k$ is odd
and $\z(k)=0$ if $r_k$ is even. We have $r_k+\z(k)\in2\NN$ for any $k$ and $r_k+\z(k)\ge r_{k+1}+\z(k+1)$ for 
$k\in[1,\t-1]$.)

Let $A_\nn$ be the set of all pairs $(r_*,p_*)\in\car\T\ti\cp$ such that $|r_*|+|p_*|=\nn$. 
We define $\io:A_\nn@>>>\cq_\nn$ by $(r_*,p_*)\m c_*$ where the multiset of entries of $c_*$ is the union of
the multiset of entries of $r_*$ with the multiset of entries of $p_*$.
Let $c_*=(c_1\ge c_2\ge\do\ge c_\t)\in\cq_\nn$. Let $K=\{k\in[1,\t];c_k\text{ is odd}\}$. We write the multiset 
$\{c_k;k\in K\}$ as a sequence $c^1\ge c^2\ge\do\ge c^t$. (We have necessarily $\t=\nn=t\mod2$.) We associate to 
$c_*$ an element $(r_*,p_*)\in\cq\T\ti\cp$ by specifying the number of times $M_e$ (resp. $N_e$) that an integer 
$e\ge1$ appears in $r_*$ (resp. $p_*$). Let $Q_e$ be the number of times that $e$ appears as an entry of $c_*$.

(i) If $e\in2\NN+1$ and $Q_e=2g+1$ then $M_e=1$, $N_e=2g$.

(ii) If $e\in2\NN+1$ and $Q_e=2g$, so that $c^d=c^{d+1}=\do=c^{d+2g-1}=e$ with $d$ even, then $M_e=2$, $N_e=2g-2$
(if $g>0$) and $M_e=N_e=0$ (if $g=0$).

(iii) If $e\in2\NN+1$ and $Q_e=2g$ so that $c^d=c^{d+1}=\do=c^{d+2g-1}=e$ with $d$ odd then $M_e=0$, $N_e=2g$.

Thus the odd entries of $r_*$ are defined. We write them in a sequence $r^1\ge r^2\ge\do\ge r^s$.

(iv) If $e\in2\NN+2$, $Q_e=2g$ and if 

($*$) $r^{2v}>e>r^{2v+1}$ for some $v$, or $e>r^1$, or $r^s>e$ (with $s$ even),
\nl
then $M_e=0,N_e=2g$. 

(v) If $e\in2\NN+2$, $Q_e=2g$  and if ($*$) does not hold, then $M_e=2g,N_e=0$. 

Now $r_*\in\cq,p_*\in\tcp$ are defined and $|r_*|+|p_*|=\nn$.

Assume that $|r_*|>0$; then from (iv) we see that the largest entry of $r_*$ is odd. Assume that $|r_*|>0$ and
$\nn$ is even; then from (iv) we see that the smallest entry of $r_*$ is odd. If $u\in[1,s-1]$ and $r^u=r^{u+1}$ 
then from (i),(ii),(iii) we see that $u$ is even. If $u\in[1,s-1]$ and there is $k'\in[1,\t]$ such that 
$r^u>r_{k'}>r^{u+1}$, then $r_{k'}$ is even and $e=r_{k'}$ is as in (v) and $u$ must be odd. We see that 
$r_*\in\car$.

We see that $c_*\m(r_*,p_*)$ is a well defined map $\io':\cq_\nn@>>>A_\nn$; moreover, 
$\io\io':\cq_\nn@>>>\cq_\nn$ is the identity map. 

We preserve the notation for $c_*,r_*,p_*$ as above (so that $(r_*,p_*)\in\io\i(c_*)$) and we assume that
$(r'_*,p'_*)\in\io\i(c_*)$. We write the odd entries of $r'_*$ in a sequence 
$r'{}^1\ge r'{}^2\ge\do\ge r'{}^{s'}$.

Let $M'_e$ (resp. $N'_e$) be the number of times that an integer $e\ge1$ appears in $r'_*$ (resp. $p'_*$). Note 
that $M'_e+N'_e=M_e+N_e$.

In the setup of (i) we have $M'_e=1,N'_e=N_e$. (Indeed, $M'_e+N'_e$ is odd, $N'_e$ is even hence $M'_e$ is odd. 
Since $M'_e$ is $0,1$ or $2$ we see that it is $1$.)

In the setup of (ii) and assuming that $g>0$ we have $M'_e=2,N'_e=N_e$ or $M'_e=0,N'_e=N_e+2$. (Indeed, 
$M'_e+N'_e$ is even, $N'_e$ is even hence $M'_e$ is even. Since $M'_e$ is $0,1$ or $2$ we see that it is $0$ or 
$2$.) If $g=0$ we have $M'_e=N'_e=0$.

In the setup of (iii) we have $M'_e=0,N'_e=N_e$. (Indeed, $M'_e+N'_e$ is even, $N'_e$ is even hence $M'_e$ is
even. Since $M'_e$ is $0,1$ or $2$ we see that it is $0$ or $2$. Assume that $M'_e=2$. Then $e=r'{}^u=r'{}^{u+1}$
with $u$ even in $[1,s'-1]$. We have $c^d=c^{d+1}=\do=c^{d+2g-1}=e$ with $d$ odd. From the definitions we see
that $u=d\mod2$ and we have a contradiction. Thus, $M'_e=0$.)

Now the sequence $r'{}^1\ge r'{}^2\ge\do\ge r'{}^{s'}$ is obtained from the sequence
$r^1\ge r^2\ge\do\ge r^s$ by deleting some pairs of the form $r^{2h}=r^{2h+1}$. Hence in the setup of (iv) we 
have $r^{2v}=r'{}^{2v'}>e>r'{}^{2v'+1}=r^{2v+1}$ for some $v'$ or $e>r'{}^1$ or $r'{}^{s'}>e$ (with $s,s'$ even)
and we see that $M'_e=0$ so that $N'_e=2g=N_e$. 

In the setup of (v) we have $N'_e\ge 0$.

We see that in all cases we have $N'_e\ge N_e$. It follows that $\sum_eN'_e\ge\sum_eN_e$ (and the equality 
implies that $N'_e=N_e$ for all $e$ hence $(r'_*,p'_*)=(r_*,p_*)$). We see that 

(a) {\it for any $c_*\in\cq_\nn$ there is exactly one element $(r_*,p_*)\in\io\i(c_*)$ such that the number of
entries of $p_*$ is minimal (that element is $\io'(c_*)$).}
\nl
Let $\dot A_\nn$ be the set of all pairs $(p'_*,p_*)\in\cp^\k\T\ti\cp$ such that $2|p'_*|+1+|p_*|=\nn$. 
We have a bijection 

(b) $\dot A_\nn@>\si>>>A_\nn$, $(p'_*,p_*)\m(\Xi(p'_*),p_*)$
\nl
where $\Xi$ is as above (with $n$ replaced by $|p'_*|$).

\subhead 1.3\endsubhead
In this subsection we assume that $\nn$ is even. 
Let $\ce_\nn$ be the set of all $p_*\in\ti\cp_\nn$ such that any entry of $p_*$ is even.
We can view $\ce_\nn$ as a subset of $A_\nn$ by $p_*\m(r_*,p_*)$ where $r_*\in\car_0$ is the empty sequence.
Let $\tA_\nn=A_\nn-\ce_\nn$. Moreover, we can view $\te_\nn$ as a subset of $\cq_\nn$ by $p_*\m p*$.
Let $\ti\cq_\nn=\cq_\nn-\ce_\nn$. Note that $\io$, $\io'$ restrict to the identity map $\ce_\nn@>>>\ce_\nn$.
Let $\ti\io:\tA_\nn@>>>\ti\cq_\nn$, $\ti\io':\ti\cq_\nn@>>>\tA_\nn$ be the
restrictions of $\io,\io'$. Note that $\ti\io\ti\io':\ti\cq_\nn@>>>\ti\cq_\nn$ is the identity map. 

We can view $\ce_\nn$ as a subset of $\dot A_\nn$ by $p_*\m(p'_*,p_*)$ where $p'_*$ is the empty sequence.
Let $\dot{\tA}_\nn=\dot A_\nn-\ce_\nn$. The bijection 1.2(b) restricts to a bijection

(a) $\dot{\tA}_\nn@>\si>>>\tA_\nn$.

\subhead 1.4\endsubhead
In this subsection we assume that $p\ne2$.
Let $V$ be a $\kk$-vector space of finite dimension $\nn\ge3$. Let $\k=0$ if $\nn$ is even, $\k=1$ if $\nn$ is
odd. Let $n=(\nn-\k)/2$. Let $\e\in\{1,-1\}$.
Assume that $V$ has a fixed nonsingular bilinear form $(,):V\T V@>>>\kk$ such that $(x,y)=\e(y,x)$ for all
$x,y\in V$. Let $Is(V)$ be the group of all isometries of $(,)$ (a closed subgroup of $GL(V)$). We assume that 
$G$ is the identity component of $Is(V)$.

Assume first that $\e=-1$. We identify $\uuG=\ct_\nn$ by associating to $\g\in\uuG$ the multiset consisting of 
the sizes of the Jordan blocks of an element of $\g$. We identify (as in \cite{\WE, 1.4, 1.5}) $\WW$ with the 
group $W$
of permutations of $[1,\nn]$ commuting with the involution $i\m\nn-i+1$. We identify $\uWW$ with $A'_\nn$ by 
associating to $C\in\uWW$ (with $w\in C$) the pair $(r_*,p_*)$ where $r_*$ is the multiset consisting of the
sizes of cycles of $w$ which commute with the involution above and $p_*$ is the multiset consisting of the sizes
of the remaining cycles of $w$. Using \cite{\WE, 3.7, 1.1} we see that the map $\Ph:\uWW@>>>\uuG$ becomes the map
$\io:A'_\nn@>>>\ct_\nn$ in 1.1 and 0.2 folows from 1.1(a).

Assume next that $\e=1$. In the case where $\nn$ is even we assume that $\nn\ge8$ and we let $\uuG_0$ be the set 
of unipotent classes in $G$ which are also conjugacy classes in $Is(V)$. We identify $\uuG=\cq_\nn$ (if $\nn$ is 
odd) and $\uuG_0=\ti\cq_\nn$ (if $\nn$ is even) by 
associating to $\g$ in $\uuG$ or $\uuG_0$ the multiset consisting of the sizes of the Jordan blocks of an element
of $\g$. If $\nn$ is odd we identify (as in \cite{\WE, 1.4, 1.5}) $\WW$ with the group $W$ of permutations of 
$[1,\nn]$ commuting with the involution $i\m\nn-i+1$. If $\nn$ is even we identify (as in \cite{\WE, 1.4, 1.5})
$\WW$ with the group $W'$ of even permutations of $[1,\nn]$ commuting with the involution $i\m\nn-i+1$; in this 
case let $\uWW_0$ be the set of conjugacy classes in $WW$ which are also conjugacy classes of the group of all
permutations of $[1,\nn]$ commuting with the involution above.

We identify $\uWW=\dot A_\nn$ (if $\nn$ is odd) and $\uWW_0=\dot{\tA}_\nn$ (if $\nn$ is even) by associating to 
$C$ in $\uWW$ or $\uWW_0$ (with $w\in C$) the pair $(p'_*,p_*)$ where $p'_*$ is the multiset consisting of the 
half sizes of cycles of $w$ (other than fixed points) which commute with the involution above and $p_*$ is the 
multiset consisting of the sizes of cycles of $w$ which do not commute with the involution above. Using 1.2(b) 
and 1.3(a) we identify $A_\nn=\dot A_\nn$ if $\nn$ is odd and $\tA_\nn=\dot{\tA}_\nn$ if $\nn$ is even. Using 
\cite{\WE, 3.8, 3.9, 1.1} we see that the map $\Ph:\uWW@>>>\uuG$ becomes the map $\io:A_\nn@>>>\cq_\nn$ in 1.2 
(if $\nn$ is odd) and the map $\Ph:\uWW_0@>>>\uuG_0$ becomes the map $\ti\io:\tA_\nn@>>>\cq_\nn$ in 1.3 (if $\nn$
is even) and 0.2 folows from 1.2(a). (Note that if $\nn$ is even and $\g\in\uuG-\uuG_0$ then $\Ph\i(\g)$ is a 
single element hence for such $\g$ the statement of 0.2 is trivial.)

\head 2. Exceptional groups\endhead
\subhead 2.1\endsubhead
In 2.2-2.6 we describe explicitly the map $\Ph:\uWW@>>>\uuG$ in the case where $G$ is a simple exceptional group 
in the form of tables. Each table consists of lines of the form $[a,b,\do,r]\m s$
where $s\in\uuG$ is specified by its name in \cite{\SPA} and $a,b,\do,r$ are the elements of $\uWW$ which are
mapped by $\Ph$ to $s$ (here $a,b,\do,r$ are specified by their name in \cite{\CA}); by inspection we see that
0.2 holds in each case and in fact $\Ps(s)=a$ is the first element of $\uWW$ in the list $a,b,\do,r$. The tables 
are obtained from the results in \cite{\WE}.

\subhead 2.2. Type $G_2$\endsubhead
If $p\ne3$ we have

 $[A_0]\m A_0$

$[A_1]\m A_1$ 

$[A_1+\tA_1,\tA_1] \m \tA_1$  

$[A_2]\m G_2(a_1)$

$[G_2]\m G_2$ 
\nl
When $p=3$ the line $[A_1+\tA_1,\tA_1] \m\tA_1$ should be replaced by $[A_1+\tA_1]\m\tA_1$, $[\tA_1]\m(\tA_1)_3$. 

\subhead 2.3. Type $F_4$\endsubhead
If $p\ne2$ we have 

$[A_0]\m A_0$

$[A_1]\m A_1$

$[2A_1,\tA_1]\m \tA_1$   

$[4A_1,3A_1,2A_1+\tA_1,A_1+\tA_1]\m A_1+\tA_1$ 

$[A_2]\m A_2$

$[\tA_2]\m\tA_2$

$[A_2+\tA_1]\m A_2+\tA_1$

$[A_2+\tA_2,\tA_2+A_1]\m \tA_2+A_1$ 

$[A_3,B_2]\m B_2$  

$[A_3+\tA_1,B_2+A_1]\m C_3(a_1)$  

$[D_4(a_1)]\m F_4(a_3)$

$[D_4,B_3] \m B_3$ 

$[C_3+A_1,C_3] \m C_3$

$[F_4(a_1)]\m F_4(a_2)$ 

$[B_4]\m F_4(a_1)$

$[F_4]\m F_4$
\nl
When $p=2$ the lines $[2A_1,\tA_1]\m\tA_1$, $[A_2+\tA_2,\tA_2+A_1]\m\tA_2+A_1)$, $[A_3,B_2]\m B_2$,
$[A_3+\tA_1,B_2+A_1]\m C_3(a_1)$, should be replaced by

$[2A_1]\m\tA_1$ , $[\tA_1]\m(\tA_1)_2$

$[A_2+\tA_2]\m\tA_2+A_1$, $[\tA_2+A_1]\m(\tA_2+A_1)_2$

$[A_3]\m B_2$,  $[B_2]\m(B_2)_2$

$[A_3+\tA_1]\m C_3(a_1)$, $[B_2+A_1]\m(C_3(a_1))_2$
\nl
respectively.

\subhead 2.4. Type $E_6$\endsubhead
We have

$[A_0]\m A_0$ 

$[A_1]\m A_1$

$[2A_1] \m2A_1$ 

$[4A_1,3A_1]\m3A_1$

$[A_2]\m A_2$

$[A_2+A_1]\m A_2+A_1$ 

$[2A_2]\m2A_2$

$[A_2+2A_1] \m A_2+2A_1$

$[A_3]\m A_3$

$[3A_2,2A_2+A_1]\m2A_2+A_1$ 

$[A_3+2A_1,A_3+A_1]\m A_3+A_1$ 

$[D_4(a_1)]\m D_4(a_1)$ 

$[A_4]\m A_4$ 

 $[D_4] \m D_4$

 $[A_4+A_1]\m A_4+A_1$    

$[A_5+A_1,A_5] \m A_5$

$[D_5(a_1)]\m D_5(a_1)$ 

$[E_6(a_2)]\m A_5+A_1$ 

$[D_5]\m D_5$ 

$[E_6(a_1)] \m E_6(a_1)$ 

$[E_6]\m E_6$ 

\subhead 2.5. Type $E_7$\endsubhead
If $p\ne 2$ we have  

$[A_0]\m A_0$

 $[A_1] \m A_1$ 

 $[2A_1] \m 2A_1$

$[3A_1'] \m 3A_1''$ 

 $[ 4A_1'',3A_1'']  \m 3A_1'$

 $[A_2] \m A_2$  

$[7A_1, 6A_1,5A_1, 4A_1']\m 4A_1$ 

 $[A_2+A_1] \m A_2+A_1$ 

 $[A_2+2A_1]\m A_2+2A_1$ 

 $[A_3] \m A_3$ 

 $[2A_2]\m 2A_2$
 
 $[A_2+3A_1]\m A_2+3A_1$  

$[A_3+A_1'] \m (A_3+A_1)''$ 

 $[3A_2,2A_2+A_1] \m 2A_2+A_1$ 

$[A_3+2A_1'',A_3+A_1'']  \m (A_3+A_1)'$ 

  $[D_4(a_1)]\m D_4(a_1)$ 

 $[A_3+3A_1,A_3+2A_1']  \m A_3+2A_1$ 

 $[D_4] \m D_4$ 

 $[D_4(a_1)+A_1]\m D_4(a_1)+A_1$ 

  $[D_4(a_1)+2A_1,A_3+A_2]  \m A_3+A_2$ 

 $[2A_3+A_1,A_3+A_2+A_1]  \m A_3+A_2+A_1$ 

$[A_4]\m A_4$ 

 $[D_4+3A_1,D_4+2A_1,D_4+A_1] \m D_4+A_1$ 

 $[A_5'] \m A_5''$ 

 $[A_4+A_1] \m A_4+A_1$ 

 $[D_5(a_1)]\m D_5(a_1)$  

 $[A_4+A_2] \m A_4+A_2 $

 $[A_5+A_1'',A_5''] \m A_5'$  

$[A_5+A_2,A_5+A_1']  \m (A_5+A_1)''$  

$[D_5(a_1)+A_1] \m D_5(a_1)+A_1$ 

$[E_6(a_2)] \m (A_5+A_1)'$ 

$[D_6(a_2)+A_1,D_6(a_2)] \m  D_6(a_2) $ 

$[E_7(a_4)]\m D_6(a_2)+A_1$ 

$[D_5]\m D_5$   

$[A_6] \m A_6$

$[D_5+A_1] \m D_5+A_1$ 

$[D_6(a_1)] \m D_6(a_1)$ 

 $[A_7] \m D_6(a_1)+A_1$ 

 $[E_6(a_1)] \m E_6(a_1)$  

 $[D_6+A_1,D_6]  \m D_6$ 

  $[E_6] \m E_6$ 

 $[E_7(a_3)] \m  D_6+A_1$  

 $[E_7(a_2)] \m E_7(a_2)$ 

 $[E_7(a_1)] \m E_7(a_1)$ 

  $[E_7] \m E_7$ 
\nl
If $p=2$, the line $[D_4(a_1)+2A_1,A_3+A_2]\m A_3+A_2)$ should be replaced by 
$[D_4(a_1)+2A_1]\m A_3+A_2$, $[A_3+A_2]\m (A_3+A_2)_2$.

\subhead 2.6. Type $E_8$ \endsubhead
If $p\ne2,3$ we have

$ [A_0] \m A_0$ 

$ [A_1] \m A_1$ 

$[2A_1] \m 2A_1$

$[4A_1', 3A_1]\m 3A_1 $

$[A_2] \m A_2$

$[8A_1,7A_1,6A_1,5A_1,4A_1''] \m 4A_1$

$[A_2+A_1]\m A_2+A_1$

$[A_2+2A_1] \m A_2+2A_1$

$[A_3]\m A_3$

$[A_2+4A_1, A_2+3A_1]\m A_2+3A_1$

$[2A_2] \m 2A_2$

$[3A_2, 2A_2+A_1]\m 2A_2+A_1$

$ [A_3+2A_1', A_3+A_1] \m A_3+A_1$

$[D_4(a_1)]\m D_4(a_1)$

$[4A_2,3A_2+A_1,2A_2+2A_1]\m 2A_2+2A_1$

$[D_4]\m D_4$

$[A_3+4A_1,A_3+3A_1,A_3+2A_1'']\m A_3+2A_1$

$[D_4(a_1)+A_1]\m D_4(a_1)+A_1$

$[2A_3',A_3+A_2\m A_3+A_2]$  

$[A_4] \m A_4$

$[2A_3+2A_1, A_3+A_2+2A_1,2A_3+A_1,A_3+A_2+A_1]\m A_3+A_2+A_1$

$[D_4(a_1)+A_2]\m D_4(a_1)+A_2$

$[D_4+4A_1, D_4+3A_1, D_4+2A_1, D_4+A_1]\m D_4+A_1$

$[2D_4(a_1), D_4(a_1)+A_3, 2A_3''] \m 2A_3$

$[A_4+A_1]\m A_4+A_1$

$[D_5(a_1)]\m D_5(a_1)$

$[A_4+2A_1] \m A_4+2A_1$

$[A_4+A_2]\m A_4+A_2$

$[A_4+A_2+A_1]\m A_4+A_2+A_1$

$[D_5(a_1)+A_1]\m D_5(a_1)+A_1$

$[A_5+A_1', A_5]\m A_5$

$[D_4+A_3,D_4+A_2]\m D_4+A_2$ 

$[E_6(a_2)]\m (A_5+A_1)''$

$[2A_4,A_4+A_3] \m A_4+A_3$

$[D_5]\m D_5$

$[D_5(a_1)+A_3, D_5(a_1)+A_2]\m D_5(a_1)+A_2$

$[A_5+A_2+A_1,A_5+A_2,A_5+2A_1,A_5+A_1''] \m (A_5+A_1)'$

$[E_6(a_2)+A_2, E_6(a_2)+A_1] \m A_5+2A_1$

$[2D_4,D_6(a_2)+A_1,D_6(a_2)] \m D_6(a_2)$

$[E_7(a_4)+A_1, E_7(a_4)] \m A_5+A_2$      

$ [D_5+2A_1,D_5+A_1] \m D_5+A_1$

$[E_8(a_8)]\m2A_4$

$[D_6(a_1)\m D_6(a_1)]$

$[A_6]\m A_6$

$[A_6+A_1]\m A_6+A_1$

$[A_7'] \m D_6(a_1)+A_1$

$[A_7+A_1,D_5+A_2]\m D_5+A_2$ 

$ [E_6(a_1)]\m E_6(a_1)$

$[D_6+2A_1, D_6+A_1 ,D_6] \m D_6$

$[D_7(a_2)]\m D_7(a_2)$

$[E_6] \m E_6$

$[D_8(a_3), A_7''] \m A_7$                      

$[ E_6(a_1)+A_1]\m E_6(a_1)+A_1 $

$[E_7(a_3)]\m D_6+A_1$

$[A_8]\m D_8(a_3)$

$[D_8(a_2), D_7(a_1)] \m D_7(a_1)$ 

$[E_6+A_2, E_6+A_1]\m E_6+A_1$

$[E_7(a_2)+A_1, E_7(a_2)]\m E_7(a_2)$

$[E_8(a_6)]\m A_8$

$[E_8(a_7)]\m E_7(a_2)+A_1$

$[D_8(a_1),D_7]\m D_7$

$[E_8(a_3)]\m D_8(a_1)$

$[E_7(a_1)]\m E_7(a_1)$

$[D_8]\m E_7(a_1)+A_1$

$[E_8(a_5)] \m D_8"$

$ [E_7+A_1, E_7] \m E_7$

$ [E_8(a_4)]  \m E_7+A_1"$

$[E_8(a_2)] \m E_8(a_2)$

$[E_8(a_1)] \m E_8(a_1)$

$[ E_8]\m E_8$
\nl
If $p=3$ the line $[D_8(a_3),A_7'']\m A_7$ should be replaced by $[D_8(a_3)] \m A_7$, $[A_7''] \m(A_7)_3 $.
If $p=2$ the lines $[2A_3',A_3+A_2]\m A_3+A_2$,  $[D_4+A_3,D_4+A_2]\m D_4+A_2$, 
$[A_7+A_1, D_5+A_2]\m D_5+A_2$,  $[D_8(a_2),D_7(a_1)] \m D_7(a_1)$ should be replaced by    

$[2A_3'] \m A_3+A_2$, $[A_3+A_2]\m(A_3+A_2)_2$

$[D_4+A_3] \m D_4+A_2$,  $[D_4+A_2]\m (D_4+A_2)_2$
 
$[A_7+A_1]\m D_5+A_2$, $[D_5+A_2]\m (D_5+A_2)_2$
   
$[D_8(a_2)]\m D_7(a_1)$, $[D_7(a_1)]\m (D_7(a_1))_2$    
\nl
respectively.

\head 3. Relation with unipotent pieces\endhead
\subhead 3.1\endsubhead
Let $G'$ be a connected reductive group over $\CC$ of the same type as $G$. We identify $\WW$ with the Weyl group
of $G'$. Let $\uuG'$ be the set of unipotent classes in $G'$. In \cite{\UNI, 6.8} we have defined a partition of
the unipotent variety of $G$ into "unipotent pieces" indexed by $\uuG'$. We define a surjective map
$\r:\uuG@>>>\uuG'$ by associating to $\g\in\uuG$ the element $\g'\in\uuG'$ which indexes the unipotent piece that
contains $\g$. Let $\p:\uuG'@>>>\uuG$ be the map defined in \cite{\WE, 4.1}. Note that $\r\p=1$.
Let $\Ph':\uWW@>>>\uuG'$, $\Ps':\uuG'@>>>uWW$ be the maps analogous to $\Ph:\uWW@>>>\uuG$ in \cite{\WE, 4.5} and 
$\Ps:\uuG@>>>\uWW$ in 0.2. One can show that

(a)   $\Ph'=\r\Ph$
\nl
and (assuming that $\Ps$ is defined):

(b) $\Ps'=\Ps\p$.
\nl
In the case where $G$ is simple of exceptional type this follows from the tables in Section 2.

\enddocument